\documentclass[thmsa,12pt]{article}%
\usepackage{amsmath}
\usepackage{amssymb}
\usepackage{sw20lart}
\usepackage{amsfonts}
\usepackage{graphicx}%
\begin{document}

\title{Metric Compatible Covariant Derivatives}
\author{{\footnotesize W. A. Rodrigues Jr}.$^{1}$,{\footnotesize \ V. V.
Fern\'{a}ndez}$^{1}${\footnotesize and A. M. Moya}$^{2}$\\$^{1}\hspace{-0.1cm}${\footnotesize Institute of Mathematics, Statistics and
Scientific Computation}\\{\footnotesize \ IMECC-UNICAMP CP 6065}\\{\footnotesize \ 13083-859 Campinas, SP, Brazil }\\{\footnotesize e-mail:walrod@ime.unicamp.br }\\{\footnotesize \ }$^{2}${\footnotesize Department of Mathematics, University
of Antofagasta, Antofagasta, Chile} \\{\footnotesize e-mail: mmoya@uantof.cl}}
\maketitle

\begin{abstract}
This paper, sixth in a series of eight, uses the geometric calculus on
manifolds developed in the previous papers of the series to introduce through
the concept of a metric extensor field $g$ a metric structure for a smooth
manifold $M$. The associated metric compatible connection extensor field, the
associated Christoffel operators \ and a notable decomposition of those
objects are given The paper introduces also the concept of a geometrical
structure for a manifold $M$ as a triple $(M,g,\gamma)$, where $\gamma$ is a
connection extensor field defining a parallelism structure for $M$. Next, the
theory of metric compatible covariant derivatives is given and a relationship
between the connection extensor fields and covariant derivatives of two
deformed (metric compatible) geometrical structures $(M,g,\gamma)$ and
$(M,\eta,\gamma^{\prime})$ is determined.

\end{abstract}
\tableofcontents

\section{Introduction}

This is the sixth paper in a series of eight. Here, using the geometric
calculus on manifolds developed in previous papers \cite{1,2} we introduce in
Section 2 a metric structure on a smooth manifold through the concept of a
metric extensor field $g$. Christoffel operators and the associated
Levi-Civita connection field are given. In Section 3 we introduce the notion
of a geometrical structure for a manifold $M$ as a triple $(M,g,\gamma)$,
where $\gamma$ is a connection extensor field defining a parallelism structure
for $M$ and study the theory of metric compatible covariant derivatives.
There, the relationship between the connection extensor fields and covariant
derivatives corresponding to \textit{deformed} (metric compatible) geometrical
structures are given. The crucial result is Theorem 2. In section 4 we present
our conclusions.

\section{Metric Structure}

Let$\ U$ be an open subset\footnote{In this paper we will use the nomenclature
and notations just used in \cite{1},\cite{2}.} of $U_{o}.$ Any symmetric and
non-degenerate smooth $(1,1)$-extensor field on $U,$ namely $g,$ will be said
to be a \emph{metric field} on $U.$ This means that $g:U_{o}\rightarrow
ext_{1}^{1}(\mathcal{U}_{o})$ satisfies $g_{(p)}=g_{(p)}^{\dagger}$ and
$\det[g]\neq0,$ for each $p\in U,$ and for all $v\in\mathcal{V}(U)$ the vector
field defined by $U_{o}\ni p\mapsto g_{(p)}(v(p))$ belongs to $\mathcal{V}%
(U),$ see \cite{1}.

The open set $U$ equipped with such a metric field $g,$ namely $(U,g),$ will
be said to be a \emph{metric structure} on $U.$

The existence of a metric field on $U$ makes possible the introduction of
three kinds of \emph{metric products} of smooth multivector fields on $U$.
These are: (a) the $g$-\emph{scalar product }of $X,Y\in\mathcal{M}(U),$ namely
$X\underset{g}{\cdot}Y\in\mathcal{S}(U)$; (b) the \emph{left} and \emph{right}
$g$-\emph{contracted products }of $X,Y\in\mathcal{M}(U),$ namely
$X\underset{g}{\lrcorner}Y\in\mathcal{M}(U)$ and $X\underset{g}{\llcorner}%
Y\in\mathcal{M}(U).$ These products are defined by
\begin{align}
(X\underset{g}{\cdot}Y)(p)  &  =\underline{g}_{(p)}(X(p))\cdot
Y(p)\label{MS.1a}\\
(X\underset{g}{\lrcorner}Y)(p)  &  =\underline{g}_{(p)}(X(p))\lrcorner
Y(p)\label{MS.1b}\\
(X\underset{g}{\llcorner}Y)(p)  &  =X(p)\llcorner\underline{g}_{(p)}%
(Y(p)),\text{ for each }p\in U. \label{MS.1c}%
\end{align}
Note that in the above formulas $\underline{g}$ is the \emph{\ extended}
(extensor field) of $g$ \cite{3}.

The $g$-\emph{Clifford product }of $X,Y\in\mathcal{M}(U),$ namely
$X\underset{g}{}Y\in\mathcal{M}(U),$ is defined by the following axioms.

For all $f\in\mathcal{S}(U),$ $b\in\mathcal{V}(U)$ and $X,Y,Z\in
\mathcal{M}(U)$
\begin{align}
f\underset{g}{}X  &  =X\underset{g}{}f=fX\text{ (scalar multiplication on
}\mathcal{M}(U)\text{).}\label{MS.2a}\\
b\underset{g}{}X  &  =b\underset{g}{\lrcorner}X+b\wedge X,\label{MS.2b}\\
X\underset{g}{}b  &  =X\underset{g}{\llcorner}b+X\wedge b.\label{MS.2c}\\
X\underset{g}{}(Y\underset{g}{}Z)  &  =(X\underset{g}{}Y)\underset{g}{}Z.
\label{MS.2d}%
\end{align}

$\mathcal{M}(U)$ equipped with each one of the products $(\underset
{g}{\lrcorner})$ or $(\underset{g}{\llcorner})$ is a non-associative algebra
induced by the respective $b$-interior algebra of multivectors. They are
called the $g$-\emph{interior algebras of smooth multivector fields.}

$\mathcal{M}(U)$ equipped with $(\underset{g}{})$ is an associative algebra
fundamentally induced by the $b$-Clifford algebra of multivectors. It is
called the $g$-\emph{Clifford algebra of smooth multivector fields.}

\subsection{Christoffel Operators}

Given a metric structure $(U,g)$ we can introduce the following two remarkable
operators which map $3$-uples of smooth vector fields into smooth scalar fields.

(a) The mapping $\left[  \left.  {}\right.  ,\left.  {}\right.  ,\left.
{}\right.  \right]  :\mathcal{V}(U)\times\mathcal{V}(U)\times\mathcal{V}%
(U)\rightarrow\mathcal{S}(U)$ defined by
\begin{align}
\left[  a,b,c\right]   &  =\frac{1}{2}(a\cdot\partial_{o}(b\underset{g}{\cdot
}c)+b\cdot\partial_{o}(c\underset{g}{\cdot}a)-c\cdot\partial_{o}(a\underset
{g}{\cdot}b)\nonumber\\
&  +c\underset{g}{\cdot}[a,b]+b\underset{g}{\cdot}[c,a]-a\underset{g}{\cdot
}[b,c]) \label{CHO.1}%
\end{align}
which is called the \emph{Christoffel operator of first kind.} Of course, it
is associated to $(U,g).$

(b) The mapping $%
\genfrac{\{}{\}}{0pt}{}{\left.  {}\right.  }{\left.  {}\right.  ,\left.
{}\right.  }%
:\mathcal{V}(U)\times\mathcal{V}(U)\times\mathcal{V}(U)\rightarrow
\mathcal{S}(U)$ defined by
\begin{equation}%
\genfrac{\{}{\}}{0pt}{0}{c}{a,b}%
=[a,b,g^{-1}(c)] \label{CHO.2}%
\end{equation}
which is called the \emph{Christoffel operator of second kind.}

We summarize here some of the basic properties of the Christoffel operator of
first kind.\vspace{0.1in}

\textbf{i.} For all $f\in\mathcal{S}(U),$ and $a,a^{\prime},b,b^{\prime
},c,c^{\prime}\in\mathcal{V}(U)$,
\begin{align}
\lbrack a+a^{\prime},b,c]  &  =[a,b,c]+[a^{\prime},b,c],\label{CHO.3a}\\
\lbrack fa,b,c]  &  =f[a,b,c].\label{CHO.3b}\\
\lbrack a,b+b^{\prime},c]  &  =[a,b,c]+[a,b^{\prime},c],\label{CHO.3c}\\
\lbrack a,fb,c]  &  =f[a,b,c]+(a\cdot\partial_{o}f)b\underset{g}{\cdot
}c.\label{CHO.3d}\\
\lbrack a,b,c+c^{\prime}]  &  =[a,b,c]+[a,b,c^{\prime}],\label{CHO.3e}\\
\lbrack a,b,fc]  &  =f[a,b,c]. \label{CHO.3f}%
\end{align}
These equations say that the Christoffel operator of the first kind has the
linearity property which respect to the first and third smooth vector field variables.

\textbf{ii.} For all $a,b,c\in\mathcal{V}(U)$,
\begin{align}
\lbrack a,b,c]+[b,a,c]  &  =a\cdot\partial_{o}(b\underset{g}{\cdot}%
c)+b\cdot\partial_{o}(c\underset{g}{\cdot}a)-c\cdot\partial_{o}(a\underset
{g}{\cdot}b)\nonumber\\
&  +b\underset{g}{\cdot}[c,a]-a\underset{g}{\cdot}[b,c],\label{CHO.4a}\\
\lbrack a,b,c]-[b,a,c]  &  =c\underset{g}{\cdot}[a,b].\label{CHO.4b}\\
\lbrack a,b,c]+[a,c,b]  &  =a\cdot\partial_{o}(b\underset{g}{\cdot
}c),\label{CHO.4c}\\
\lbrack a,b,c]-[a,c,b]  &  =b\cdot\partial_{o}(c\underset{g}{\cdot}%
a)-c\cdot\partial_{o}(a\underset{g}{\cdot}b)\nonumber\\
&  +c\underset{g}{\cdot}[a,b]+b\underset{g}{\cdot}[c,a]-a\underset{g}{\cdot
}[b,c].\label{CHO.4d}\\
\lbrack a,b,c]+[c,b,a]  &  =b\cdot\partial_{o}(c\underset{g}{\cdot
}a)+c\underset{g}{\cdot}[a,b]-a\underset{g}{\cdot}[b,c],\label{CHO.4e}\\
\lbrack a,b,c]-[c,b,a]  &  =a\cdot\partial_{o}(b\underset{g}{\cdot}%
c)-c\cdot\partial_{o}(a\underset{g}{\cdot}b)+b\underset{g}{\cdot}[c,a].
\label{CHO.4f}%
\end{align}

\subsection{Levi-Civita Connection Field}

\textbf{Theorem}. There exists a smooth $(1,2)$-extensor field on $U,$ namely
$\omega_{0},$ such that the Christoffel operator of first kind can be written
as
\begin{equation}
\lbrack a,b,c]=(a\cdot\partial_{o}b+\frac{1}{2}g^{-1}\circ(a\cdot\partial
_{o}g)(b)+\omega_{0}(a)\underset{g}{\times}b)\underset{g}{\cdot}c.
\label{LCC.1}%
\end{equation}
Such $\omega_{0}$ is given by
\begin{equation}
\omega_{0}(a)=-\frac{1}{4}\underline{g}^{-1}(\partial_{b}\wedge\partial
_{c})a\cdot((b\cdot\partial_{o}g)(c)-(c\cdot\partial_{o}g)(b)). \label{LCC.2}%
\end{equation}

\textbf{Proof}

By using $a\cdot\partial_{o}(X\underset{g}{\cdot}Y)=(a\cdot\partial
_{o}X)\underset{g}{\cdot}Y+X\underset{g}{\cdot}(a\cdot\partial_{o}%
Y)+(a\cdot\partial_{o}\underline{g})(X)\cdot Y,$ for all $X,Y\in
\mathcal{M}(U),$ we have
\begin{equation}
a\cdot\partial_{o}(b\underset{g}{\cdot}c)=(a\cdot\partial_{o}b)\underset
{g}{\cdot}c+b\underset{g}{\cdot}(a\cdot\partial_{o}c)+(a\cdot\partial
_{o}g)(b)\cdot c, \label{LCC.N1}%
\end{equation}
and, by cycling the letters $a,b$ and $c$, we get
\begin{align}
b\cdot\partial_{o}(c\underset{g}{\cdot}a)  &  =(b\cdot\partial_{o}%
c)\underset{g}{\cdot}a+c\underset{g}{\cdot}(b\cdot\partial_{o}a)+(b\cdot
\partial_{o}g)(c)\cdot a,\label{LCC.N2}\\
-c\cdot\partial_{o}(a\underset{g}{\cdot}b)  &  =-(c\cdot\partial
_{o}a)\underset{g}{\cdot}b-a\underset{g}{\cdot}(c\cdot\partial_{o}%
b)-(c\cdot\partial_{o}g)(a)\cdot b. \label{LCC.N3}%
\end{align}

A straightforward calculation yields
\begin{align}
c\underset{g}{\cdot}[a,b]  &  =c\underset{g}{\cdot}(a\cdot\partial
_{o}b)-c\underset{g}{\cdot}(b\cdot\partial_{o}a),\label{LCC.N4}\\
b\underset{g}{\cdot}[c,a]  &  =b\underset{g}{\cdot}(c\cdot\partial
_{o}a)-b\underset{g}{\cdot}(a\cdot\partial_{o}c),\label{LCC.N5}\\
-a\underset{g}{\cdot}[b,c]  &  =-a\underset{g}{\cdot}(b\cdot\partial
_{o}c)+a\underset{g}{\cdot}(c\cdot\partial_{o}b). \label{LCC.N6}%
\end{align}

Now, by adding Eqs.(\ref{LCC.N1}), (\ref{LCC.N2}), (\ref{LCC.N3}) and
Eqs.(\ref{LCC.N4}), (\ref{LCC.N5}), (\ref{LCC.N6}) we get
\[
2[a,b,c]=2(a\cdot\partial_{o}b)\underset{g}{\cdot}c+(a\cdot\partial
_{o}g)(b)\cdot c+(b\cdot\partial_{o}g)(c)\cdot a-(c\cdot\partial_{o}g)(a)\cdot
b,
\]
hence, by taking into account the symmetry property $(n\cdot\partial
_{o}g)^{\dagger}=n\cdot\partial_{o}g,$ it follows
\begin{equation}
\lbrack a,b,c]=(a\cdot\partial_{o}b)\underset{g}{\cdot}c+\frac{1}{2}%
g^{-1}\circ(a\cdot\partial_{o}g)(b)\underset{g}{\cdot}c+\frac{1}{2}%
a\cdot((b\cdot\partial_{o}g)(c)-(c\cdot\partial_{o}g)(b)). \label{LCC.N7}%
\end{equation}

On another side, a straightforward calculation yields
\begin{align*}
\omega_{0}(a)\underset{g}{\times}b  &  =-g(b)\lrcorner\omega_{0}(a)\\
&  =\frac{1}{4}g(b)\lrcorner\underline{g}^{-1}(\partial_{p}\wedge\partial
_{q})a\cdot((p\cdot\partial_{o}g)(q)-(q\cdot\partial_{o}g)(p))\\
&  =\frac{1}{4}(b\cdot\partial_{p}g^{-1}(\partial_{q})-b\cdot\partial
_{q}g^{-1}(\partial_{p}))a\cdot((p\cdot\partial_{o}g)(q)-(q\cdot\partial
_{o}g)(p))\\
&  =\frac{1}{2}b\cdot\partial_{p}g^{-1}(\partial_{q})a\cdot((p\cdot
\partial_{o}g)(q)-(q\cdot\partial_{o}g)(p))\\
&  =\frac{1}{2}g^{-1}(\partial_{q})a\cdot((b\cdot\partial_{o}g)(q)-(q\cdot
\partial_{o}g)(b)),
\end{align*}
hence, it follows
\begin{equation}
(\omega_{0}(a)\underset{g}{\times}b)\underset{g}{\cdot}c=\frac{1}{2}%
a\cdot((b\cdot\partial_{o}g)(c)-(c\cdot\partial_{o}g)(b)). \label{LCC.N8}%
\end{equation}

Finally, putting Eq.(\ref{LCC.N8}) into Eq.(\ref{LCC.N7}), we get the required
result.$\blacksquare$

The \emph{smooth vector elementary} $2$-\emph{extensor field} on $U,$ namely
$\lambda,$ defined by
\begin{equation}
\lambda(a,b)=\frac{1}{2}g^{-1}\circ(a\cdot\partial_{o}g)(b)+\omega
_{0}(a)\underset{g}{\times}b \label{LCC.3}%
\end{equation}
is a well-defined connection field on $U.$ It will be called the
\emph{Levi-Civita connection field} on $U.$ The open set $U$ endowed with
$\lambda,$ namely $(U,\lambda),$ will be said to be the \emph{Levi-Civita
parallelism structure} on $U.$

The $a$-\emph{DCDO's} associated to $(U,\lambda),$ namely $D_{a}^{+}$ and
$D_{a}^{-},$ are said to be \emph{Levi-Civita }$a$-\emph{DCDO's.} They are
fundamentally defined by $D_{a}^{\pm}:\mathcal{M}(U)\rightarrow\mathcal{M}(U)$
such that
\begin{align}
D_{a}^{+}X  &  =a\cdot\partial_{o}X+\Lambda_{a}(X),\label{LCC.3a1}\\
D_{a}^{-}X  &  =a\cdot\partial_{o}X-\Lambda_{a}^{\dagger}(X), \label{LCC.3a2}%
\end{align}
where $\Lambda_{a}$ is the \emph{generalized}\footnote{We recall from \cite{3}
that the generalized of $t\in ext_{1}^{1}(\mathcal{U}_{o})$ is $T\in
ext(\mathcal{U}_{o})$ given by $T(X)=t(\partial_{n})\wedge(n\lrcorner X).$}
(extensor field) of $\lambda_{a}.$

It should be noted that such a $a$-\emph{DCDO} $D_{a}^{+}$ satisfies the
fundamental property
\begin{equation}
(D_{a}^{+}b)\cdot c=%
\genfrac{\{}{\}}{0pt}{}{c}{a,b}%
,\text{ for all }a,b,c\in\mathcal{V}(U). \label{LCC.3a3}%
\end{equation}
Eq.(\ref{LCC.3a3}) follows immediately from Eq.(\ref{LCC.1}) once we change
$c$ for $g^{-1}(c)$ and take into account the definitions given by
Eq.(\ref{CHO.2}), Eq.(\ref{LCC.3}) and Eq.(\ref{LCC.3a1}).

We present now two remarkable properties of $\omega_{0}.\vspace{0.1in}$

\textbf{i.} For all $a,b,c\in\mathcal{V}(U)$ we have the \emph{cyclic
property}
\begin{equation}
\omega_{0}(a)\underset{g}{\times}b\underset{g}{\cdot}c+\omega_{0}%
(b)\underset{g}{\times}c\underset{g}{\cdot}a+\omega_{0}(c)\underset{g}{\times
}a\underset{g}{\cdot}b=0. \label{LCC.3a}%
\end{equation}

\textbf{Proof}

By recalling Eq.(\ref{LCC.N8}) used in the proof of Eq.(\ref{LCC.1}) we can
write
\begin{equation}
\omega_{0}(a)\underset{g}{\times}b\underset{g}{\cdot}c=\frac{1}{2}%
a\cdot((b\cdot\partial_{o}g)(c)-(c\cdot\partial_{o}g)(b)), \label{LCC.M1}%
\end{equation}
and, by cycling the letters $a,b$ and $c$, we get
\begin{align}
\omega_{0}(b)\underset{g}{\times}c\underset{g}{\cdot}a  &  =\frac{1}{2}%
b\cdot((c\cdot\partial_{o}g)(a)-(a\cdot\partial_{o}g)(c)),\label{LCC.M2}\\
\omega_{0}(c)\underset{g}{\times}a\underset{g}{\cdot}b  &  =\frac{1}{2}%
c\cdot((a\cdot\partial_{o}g)(b)-(b\cdot\partial_{o}g)(a)). \label{LCC.M3}%
\end{align}

Now, by adding Eqs.(\ref{LCC.M1}), (\ref{LCC.M2}) and (\ref{LCC.M3}) we get
\begin{align*}
&  \omega_{0}(a)\underset{g}{\times}b\underset{g}{\cdot}c+\omega
_{0}(b)\underset{g}{\times}c\underset{g}{\cdot}a+\omega_{0}(c)\underset
{g}{\times}a\underset{g}{\cdot}b\\
&  =\frac{1}{2}(a\cdot(b\cdot\partial_{o}g)(c)-c\cdot(b\cdot\partial
_{o}g)(a))+\frac{1}{2}(b\cdot(c\cdot\partial_{o}g)(a)-a\cdot(c\cdot
\partial_{o}g)(b))\\
&  +\frac{1}{2}(c\cdot(a\cdot\partial_{o}g)(b)-b\cdot(a\cdot\partial
_{o}g)(c)).
\end{align*}
Then, by taking into account the symmetry property $(n\cdot\partial
_{o}g)^{\dagger}=n\cdot\partial_{o}g,$ the expected result immediately
follows.$\blacksquare$

\textbf{ii.} $\omega_{0}$ is just the $g$-\emph{gauge connection field}
associated to $(U,\lambda),$ i.e.,
\begin{equation}
\omega_{0}(a)=\frac{1}{2}\underset{g}{biv}[\lambda_{a}]. \label{LCC.3b}%
\end{equation}

\textbf{Proof}

By using the noticeable formulas: $(a\cdot\partial_{o}\tau)\circ\tau^{-1}%
+\tau\circ(a\cdot\partial_{o}\tau^{-1})=0,$ for all \emph{non-singular} smooth
$(1,1)$-extensor field $\tau,$ $(a\cdot\partial_{o}g^{-1})^{\dagger}%
=a\cdot\partial_{o}g^{-1},$ $\partial_{n}\wedge(s(n))=0,$ for all
\emph{symmetric} $s\in ext_{1}^{1}(\mathcal{U}_{o})$, and $\partial_{n}%
\wedge(B\times n)=-2B,$ where $B\in\bigwedge^{2}\mathcal{U}_{o}.$ A
straightforward calculation\footnote{Recall that [...] $biv[t]=-\partial
_{n}\wedge(t(n))$ and $\underset{g}{biv}[t]=biv[t\circ g^{-1}].$} allows us to
get
\begin{align*}
\underset{g}{biv}[\lambda_{a}]  &  =-\partial_{n}\wedge(\lambda_{a}\circ
g^{-1}(n))\\
&  =-\frac{1}{2}\partial_{n}\wedge g^{-1}\circ(a\cdot\partial_{o}g)\circ
g^{-1}(n)-\partial_{n}\wedge(\omega(a)\underset{g}{\times}g^{-1}(n))\\
&  =\frac{1}{2}\partial_{n}\wedge(a\cdot\partial_{o}g^{-1})(n)-\partial
_{n}\wedge(\omega(a)\times n),\\
&  =0+2\omega_{0}(a).\blacksquare
\end{align*}

We present now three remarkable properties of the Levi-Civita connection field
$\lambda$.\vspace{0.1in}

\textbf{i.} $\lambda$ is \emph{symmetric} with respect to the interchanging of
vector variables, i.e.,
\begin{equation}
\lambda(a,b)=\lambda(b,a). \label{LCC.3c}%
\end{equation}

\textbf{Proof}

Let us take $a,b,c\in\mathcal{V}(U),$ we have
\begin{equation}
\lambda(a,b)\underset{g}{\cdot}c=\frac{1}{2}(a\cdot\partial_{o}g)(b)\cdot
c+\omega_{0}(a)\underset{g}{\times}b\underset{g}{\cdot}c, \label{LCC.O1}%
\end{equation}
but, by interchanging the letters $a$ and $b,$ it holds
\begin{equation}
\lambda(b,a)\underset{g}{\cdot}c=\frac{1}{2}(b\cdot\partial_{o}g)(a)\cdot
c+\omega_{0}(b)\underset{g}{\times}a\underset{g}{\cdot}c. \label{LCC.O2}%
\end{equation}

Now, subtracting Eq.(\ref{LCC.O2}) from Eq.(\ref{LCC.O1}), and taking into
account Eq.(\ref{LCC.M3}) used in the proof of Eq.(\ref{LCC.3a}), we get
\[
(\lambda(a,b)-\lambda(b,a))\underset{g}{\cdot}c=\omega_{0}(c)\underset
{g}{\times}a\underset{g}{\cdot}b+\omega_{0}(a)\underset{g}{\times}%
b\underset{g}{\cdot}c-\omega_{0}(b)\underset{g}{\times}a\underset{g}{\cdot}c.
\]
Then, by recalling the multivector identity $B\underset{g}{\times}%
v\underset{g}{\cdot}w=-B\underset{g}{\times}w\underset{g}{\cdot}v,$ where
$B\in\bigwedge^{2}\mathcal{U}_{o}$ and $v,w\in\mathcal{U}_{o},$ and
Eq.(\ref{LCC.3a}), the required result immediately follows by the
non-degeneracy of the $g$-scalar product.$\blacksquare$

\textbf{ii.} The $g$-\emph{symmetric} and $g$-\emph{skew symmetric} parts of
$\lambda_{a},$ namely $\lambda_{a\pm(g)}=\dfrac{1}{2}(\lambda_{a}\pm
\lambda_{a}^{\dagger(g)}),$ are given by
\begin{align}
\lambda_{a+(g)}(b)  &  =\dfrac{1}{2}g^{-1}\circ(a\cdot\partial_{o}%
g)(b),\label{LCC.3d}\\
\lambda_{a-(g)}(b)  &  =\omega_{0}(a)\underset{g}{\times}b. \label{LCC.3f}%
\end{align}

\textbf{Proof}

We will calculate first the metric adjoint of $\lambda_{a},$ namely
$\lambda_{a}^{\dagger(g)},$ by using the fundamental property \cite{4} of the
metric adjoint operator $\left.  {}\right.  ^{\dagger(g)}$. By recalling the
symmetry property $(a\cdot\partial_{o}g)^{\dagger}=a\cdot\partial_{o}g$ and
the multivector identity $B\underset{g}{\times}v\underset{g}{\cdot
}w=-B\underset{g}{\times}w\underset{g}{\cdot}v,$ where $B\in\bigwedge
^{2}\mathcal{U}_{o}$ and $v,w\in\mathcal{U}_{o},$ we can write that
\begin{align*}
\lambda_{a}^{\dagger(g)}(b)\underset{g}{\cdot}c  &  =b\underset{g}{\cdot
}\lambda_{a}(c)\\
&  =\frac{1}{2}b\cdot(a\cdot\partial_{o}g)(c)-\omega_{0}(a)\underset{g}%
{\times}b\underset{g}{\cdot}c\\
&  =(\frac{1}{2}g^{-1}\circ(a\cdot\partial_{o}g)(b)-\omega_{0}(a)\underset
{g}{\times}b)\underset{g}{\cdot}c,
\end{align*}
hence, by the non-degeneracy of the $g$-scalar product, it follows that
\[
\lambda_{a}^{\dagger(g)}(b)=\frac{1}{2}g^{-1}\circ(a\cdot\partial
_{o}g)(b)-\omega_{0}(a)\underset{g}{\times}b.
\]
Now, we can get that
\begin{align*}
\lambda_{a}(b)+\lambda_{a}^{\dagger(g)}(b)  &  =g^{-1}\circ(a\cdot\partial
_{o}g)(b),\\
\lambda_{a}(b)-\lambda_{a}^{\dagger(g)}(b)  &  =2\omega_{0}(a)\underset
{g}{\times}b.\blacksquare
\end{align*}

\textbf{iii.} The \emph{generalized }of $\lambda_{a},$ namely $\Lambda_{a}, $
is given by the following formula
\begin{equation}
\Lambda_{a}(X)=\dfrac{1}{2}\underline{g}^{-1}\circ(a\cdot\partial
_{o}\underline{g})(X)+\omega_{0}(a)\underset{g}{\times}X, \label{LCC.4}%
\end{equation}
where $\underline{g}$ and $\underline{g}^{-1}$ are \emph{extended} of $g$ and
$g^{-1}$, respectively.

\textbf{Proof}

The above property is an immediate consequence of using the noticeable
formulas: $\tau^{-1}\circ(a\cdot\partial_{o}\tau)(\partial_{n})\wedge
(n\lrcorner X)=\underline{\tau}^{-1}\circ(a\cdot\partial_{o}\underline{\tau
})(X),$ for all \emph{non-singular} smooth $(1,1)$-extensor field $\tau,$ and
$(B\underset{g}{\times}\partial_{n})\wedge(n\lrcorner X)=B\underset{g}{\times
}X,$ where $B\in\bigwedge^{2}\mathcal{U}_{o}$ and $X\in\bigwedge
\mathcal{U}_{o}$.$\blacksquare$

\section{Geometric Structure}

A \emph{parallelism structure }$(U,\gamma)$ (see \cite{2}) is said to be
\emph{compatible} with a \emph{metric structure} $(U,g)$ if and only if
\begin{equation}
\gamma_{a+(g)}=\frac{1}{2}g^{-1}\circ(a\cdot\partial_{o}g), \label{GS.1}%
\end{equation}
i.e., $g\circ\gamma_{a}+\gamma_{a}^{\dagger}\circ g=a\cdot\partial_{o}g $.

Sometimes for abuse of language we will say that $\gamma$ is \emph{metric
compatible} (or $g$-\emph{compatible}, for short).

The open set $U$ equipped with a connection field $\gamma$ and a metric field
$g,$ namely $(U,\gamma,g),$ such that $(U,\gamma)$ is compatible with $(U,g),$
will be said to be a \emph{geometric structure} on $U.$

The Levi-Civita parallelism structure $(U,\lambda),$ according to
Eq.(\ref{LCC.3d}), is compatible with the metric structure $(U,g),$ i.e.,
$\lambda$ is $g$-compatible. It immediately follows that $(U,\lambda,g)$ is a
well-defined geometric structure on $U.$

Using the above results we have in any geometric structure $(U,\gamma,g)$ a
fundamental theorem.

\textbf{Theorem 1. }There exists a smooth $(1,2)$-extensor field on $U,$
namely $\omega,$ such that
\begin{equation}
\gamma_{a}(b)=\frac{1}{2}g^{-1}\circ(a\cdot\partial_{o}g)(b)+\omega
(a)\underset{g}{\times}b. \label{GS.1a}%
\end{equation}

We give now three properties involving $\gamma_{a}$ and $\omega.\vspace
{0.1in}$

\textbf{i.} The $g$-\emph{symmetric} and $g$-\emph{skew-symmetric parts} of
$\gamma_{a},$ namely $\gamma_{a\pm(g)}=\dfrac{1}{2}(\gamma_{a}\pm\gamma
_{a}^{\dagger}),$ are given by
\begin{align}
\gamma_{a+(g)}(b)  &  =\frac{1}{2}g^{-1}\circ(a\cdot\partial_{o}%
g)(b),\label{GS.1b}\\
\gamma_{a-(g)}(b)  &  =\omega(a)\underset{g}{\times}b. \label{GS.1c}%
\end{align}

\textbf{ii.} The \emph{generalized} of $\gamma_{a},$ namely $\Gamma_{a},$ is
given by
\begin{equation}
\Gamma_{a}(X)=\frac{1}{2}\underline{g}^{-1}\circ(a\cdot\partial_{o}%
\underline{g})(X)+\omega(a)\underset{g}{\times}X. \label{GS.1d}%
\end{equation}

\textbf{iii.} $\omega$ is just the $g$-\emph{gauge connection field}
associated to $(U,\gamma),$ i.e.,
\begin{equation}
\omega(a)=\frac{1}{2}\underset{g}{biv}[\gamma_{a}]. \label{GS.1e}%
\end{equation}

\subsection{Metric Compatible Covariant Derivatives}

The pair of $a$-\emph{DCDO's }associated to $(U,\gamma),$ namely $(\nabla
_{a}^{+},\nabla_{a}^{-}),$ is said to be \emph{metric} \emph{compatible} (or
$g$-\emph{compatible}, for short) if and only if
\begin{equation}
\nabla_{a}^{++}g=0, \label{MCD.1}%
\end{equation}
or equivalently,
\begin{equation}
\nabla_{a}^{--}g^{-1}=0. \label{MCD.1a}%
\end{equation}

We emphasize that Eq.(\ref{MCD.1}) and Eq.(\ref{MCD.1a}) are completely
equivalent to each other. It follows from the remarkable formula $(\nabla
_{a}^{++}\tau)\circ\tau^{-1}+\tau\circ(\nabla_{a}^{--}\tau^{-1})=0$, valid for
all non-singular smooth $(1,1)$-extensor field $\tau.$

$(U,\gamma)$ is compatible with $(U,g)$ if and only if $(\nabla_{a}^{+}%
,\nabla_{a}^{-})$ is metric compatible.

Indeed, let us take $b\in\mathcal{V}(U),$ we have that
\begin{align}
(\nabla_{a}^{++}g)(b)  &  =\nabla_{a}^{-}g(b)-g(\nabla_{a}^{+}b)\nonumber\\
&  =a\cdot\partial_{o}g(b)-\gamma_{a}^{\dagger}\circ g(b)-g(a\cdot\partial
_{o}b)-g\circ\gamma_{a}(b),\nonumber\\
&  =(a\cdot\partial_{o}g)(b)-g\circ\gamma_{a}(b)-\gamma_{a}^{\dagger}\circ
g(b). \label{MCD.1b}%
\end{align}
Now, if $\gamma$ is $g$-compatible, by using Eq.(\ref{GS.1}) into
Eq.(\ref{MCD.1b}), it follows that $(\nabla_{a}^{++}g)(b)=0,$ i.e.,
$(\nabla_{a}^{+},\nabla_{a}^{-})$ is $g$-compatible. And, if $(\nabla_{a}%
^{+},\nabla_{a}^{-})$ is $g$-compatible, by using Eq.(\ref{MCD.1}) into
Eq.(\ref{MCD.1b}), we get that $g\circ\gamma_{a}(b)+\gamma_{a}^{\dagger}\circ
g(b)=(a\cdot\partial_{o}g)(b),$ i.e., $\gamma$ is $g$-compatible.

We now present some basic properties which are satisfied by a $g$-compatible
pair of $a$-\emph{DCDO's,} namely $(\mathcal{D}_{a}^{+},\mathcal{D}_{a}^{-}%
)$.\vspace{0.1in}

\textbf{i.} For any $(\mathcal{D}_{a}^{+},\mathcal{D}_{a}^{-})$ we have
\begin{align}
\mathcal{D}_{a}^{++}\underline{g}  &  =0,\label{MCD.2}\\
\mathcal{D}_{a}^{\_\_}\underline{g}^{-1}  &  =0, \label{MCD.2a}%
\end{align}
where $\underline{g}$ and $\underline{g}^{-1}$are the so-called
\emph{extended} of $g$ and $g^{-1},$ respectively.

\textbf{Proof}

In order to prove the first statement we only need to check that for all
$f\in\mathcal{S}(U)$ and $b_{1},\ldots,b_{k}\in\mathcal{V}(U)$
\[
(\mathcal{D}_{a}^{++}\underline{g})(f)=0\text{ and }(\mathcal{D}_{a}%
^{++}\underline{g})(b_{1}\wedge\ldots\wedge b_{k})=0.
\]

But, by using the fundamental property $\underline{g}(f)=f,$ we get
\[
(\mathcal{D}_{a}^{++}\underline{g})(f)=\mathcal{D}_{a}^{-}\underline
{g}(f)-\underline{g}(\mathcal{D}_{a}^{+}f)=\mathcal{D}_{a}^{-}f-\underline
{g}(a\cdot\partial_{o}f)=a\cdot\partial_{o}f-a\cdot\partial_{o}f=0.
\]

And, by using the fundamental property $\underline{g}(b_{1}\wedge\ldots\wedge
b_{k})=g(b_{1})\wedge\ldots\wedge g(b_{k})$ we get
\begin{align*}
&  (\mathcal{D}_{a}^{++}\underline{g})(b_{1}\wedge\ldots\wedge b_{k})\\
&  =\mathcal{D}_{a}^{-}\underline{g}(b_{1}\wedge\ldots\wedge b_{k}%
)-\underline{g}(\mathcal{D}_{a}^{+}(b_{1}\wedge\ldots\wedge b_{k}))\\
&  =\mathcal{D}_{a}^{-}g(b_{1})\wedge\ldots g(b_{k})+\cdots+g(b_{1}%
)\wedge\ldots\mathcal{D}_{a}^{-}g(b_{k})\\
&  -g(\mathcal{D}_{a}^{+}b_{1})\wedge\ldots g(b_{k})-\ldots-g(b_{1}%
)\wedge\ldots g(\mathcal{D}_{a}^{+}b_{k})\\
&  =(\mathcal{D}_{a}^{++}g)(b_{1})\wedge\ldots g(b_{k})+\cdots+g(b_{1}%
)\wedge\ldots(\mathcal{D}_{a}^{++}g)(b_{k})=0,
\end{align*}
since $\mathcal{D}_{a}^{++}g=0.$

The second statement can be proved analogously.$\blacksquare$

\textbf{ii. }$(\mathcal{D}_{a}^{+},\mathcal{D}_{a}^{-})$ satisfies the
fundamental property
\begin{equation}
\mathcal{D}_{a}^{-}\underline{g}(X)=\underline{g}(\mathcal{D}_{a}^{+}X).
\label{MCD.3}%
\end{equation}

\textbf{Proof}

By Eq.(\ref{MCD.2}) we have that for all $X\in\mathcal{M}(U)$
\begin{align*}
(\mathcal{D}_{a}^{++}\underline{g})(X)  &  =0,\\
\mathcal{D}_{a}^{-}\underline{g}(X)-\underline{g}(\mathcal{D}_{a}^{+}X)  &
=0.\blacksquare
\end{align*}

\textbf{iii. }\emph{Ricci-like theorems}. Let $X,Y$ be a smooth multivector
fields. Then,
\begin{align}
a\cdot\partial_{o}(X\underset{g}{\cdot}Y)  &  =(\mathcal{D}_{a}^{+}%
X)\underset{g}{\cdot}Y+X\underset{g}{\cdot}(\mathcal{D}_{a}^{+}%
Y),\label{MCD.4}\\
a\cdot\partial_{o}(X\underset{g^{-1}}{\cdot}Y)  &  =(\mathcal{D}_{a}%
^{-}X)\underset{g^{-1}}{\cdot}Y+X\underset{g^{-1}}{\cdot}(\mathcal{D}_{a}%
^{-}Y). \label{MCD.4b}%
\end{align}

\textbf{Proof}

As we know, $(\mathcal{D}_{a}^{+},\mathcal{D}_{a}^{-})$ must satisfy the
fundamental property
\[
(\mathcal{D}_{a}^{+}X)\cdot Y+X\cdot(\mathcal{D}_{a}^{-}Y)=a\cdot\partial
_{o}(X\cdot Y).
\]
Then, by substituting $Y$ for $\underline{g}(Y)$ and using Eq.(\ref{MCD.3}),
the first statement follows immediately. To prove get the second statement it
is enough to substitute $Y$ for $\underline{g}^{-1}(Y)$ and once again use
Eq.(\ref{MCD.3}).$\blacksquare$

\textbf{iii. }$(\mathcal{D}_{a}^{+},\mathcal{D}_{a}^{-})$ satisfies
\emph{Leibnitz-like rules} for all of the $g$ and $g^{-1}$ suitable
products\footnote{Here, as in \cite{4} $\underset{g}{*}$ means any product,
either $(\wedge),$ $(\underset{g}{\cdot}),$ $(\underset{g}{\lrcorner
},\underset{g}{\llcorner})$ or $(g$-\emph{Clifford product}$).$ Analogously
for $\underset{g^{-1}}{*}. $}, namely $\underset{g}{*}$ and $\underset{g^{-1}%
}{*},$ of smooth multivector fields
\begin{align}
\mathcal{D}_{a}^{+}(X\underset{g}{*}Y)  &  =(\mathcal{D}_{a}^{+}X)\underset
{g}{*}Y+X\underset{g}{*}(\mathcal{D}_{a}^{+}Y),\label{MCD.5}\\
\mathcal{D}_{a}^{-}(X\underset{g^{-1}}{*}Y)  &  =(\mathcal{D}_{a}%
^{-}X)\underset{g^{-1}}{*}Y+X\underset{g^{-1}}{*}(\mathcal{D}_{a}^{-}Y).
\label{MCD.5a}%
\end{align}

\textbf{Proof}

We prove only the first statement. The other proof is analogous.

Firstly, if $\underset{g}{\ast}$ is just$(\wedge),$ then Eq.(\ref{MCD.5}) is
nothing more than the Leibnitz rule for the exterior product of smooth
multivector fields i.e.,
\begin{equation}
\mathcal{D}_{a}^{+}(X\wedge Y)=(\mathcal{D}_{a}^{+}X)\wedge Y+X\wedge
(\mathcal{D}_{a}^{+}Y). \label{MCD.5i}%
\end{equation}
As we know it is true.

Secondly, if $\underset{g}{\ast}$ is $(\underset{g}{\cdot})$, then
$\mathcal{D}_{a}^{+}(X\underset{g}{\ast}Y)=a\cdot\partial_{o}(X\underset
{g}{\cdot}Y)$ and it follows that Eq.(\ref{MCD.5}) is nothing more than the
Ricci-like theorem for $\mathcal{D}_{a}^{+}.$

In order to prove Eq.(\ref{MCD.5}) whenever $\underset{g}{\ast}$ is either
$(\underset{g}{\lrcorner})$ or $(\underset{g}{\llcorner}),$ we use the
identities $(X\underset{g}{\lrcorner}Y)\underset{g}{\cdot}Z=Y\underset
{g}{\cdot}(\widetilde{X}\wedge Z)$ and $(X\underset{g}{\llcorner}%
Y)\underset{g}{\cdot}Z=X\underset{g}{\cdot}(Z\wedge\widetilde{Y}),$ for all
$X,Y,Z\in\mathcal{M}(U),$ and Eq.(\ref{MCD.4}) and Eq.(\ref{MCD.5i}). For
instance, for the left $g$-contracted product we can indeed write
\begin{align*}
&  a\cdot\partial_{o}((X\underset{g}{\lrcorner}Y)\underset{g}{\cdot}%
Z)=a\cdot\partial_{o}(Y\underset{g}{\cdot}(\widetilde{X}\wedge Z))\\
&  \Rightarrow(\mathcal{D}_{a}^{+}(X\underset{g}{\lrcorner}Y))\underset
{g}{\cdot}Z+(X\underset{g}{\lrcorner}Y)\underset{g}{\cdot}(\mathcal{D}_{a}%
^{+}Z)\\
&  =(\mathcal{D}_{a}^{+}Y)\underset{g}{\cdot}(\widetilde{X}\wedge
Z)+Y\underset{g}{\cdot}((\mathcal{D}_{a}^{+}\widetilde{X})\wedge
Z)+Y\underset{g}{\cdot}(\widetilde{X}\wedge(\mathcal{D}_{a}^{+}Z))\\
&  \Rightarrow(\mathcal{D}_{a}^{+}(X\underset{g}{\lrcorner}Y))\underset
{g}{\cdot}Z=((\mathcal{D}_{a}^{+}X)\underset{g}{\lrcorner}Y+X\underset
{g}{\lrcorner}(\mathcal{D}_{a}^{+}Y))\underset{g}{\cdot}Z.
\end{align*}
Hence, by the non-degeneracy of the $g$-scalar product, the Leibniz rule for
$(\underset{g}{\lrcorner})$ immediately follows, i.e.,
\begin{equation}
\mathcal{D}_{a}^{+}(X\underset{g}{\lrcorner}Y)=(\mathcal{D}_{a}^{+}%
X)\underset{g}{\lrcorner}Y+X\underset{g}{\lrcorner}(\mathcal{D}_{a}^{+}Y).
\label{MCD.5ii}%
\end{equation}

In order to prove Eq.(\ref{MCD.5}) whenever $\underset{g}{\ast}$ means
$(\underset{g}{}),$ we only need to check that for all $f\in\mathcal{S}(U)$
and $b_{1},\ldots,b_{k}\in\mathcal{V}(U)$
\begin{align}
\mathcal{D}_{a}^{+}(f\underset{g}{}X)  &  =(\mathcal{D}_{a}^{+}f)\underset
{g}{}X+f\underset{g}{}(\mathcal{D}_{a}^{+}X),\label{MCD.5iii}\\
\mathcal{D}_{a}^{+}(b_{1}\underset{g}{}\ldots\underset{g}{}b_{k}\underset{g}%
{}X)  &  =(\mathcal{D}_{a}^{+}(b_{1}\underset{g}{}\ldots\underset{g}{}%
b_{k}))\underset{g}{}X\nonumber\\
&  +b_{1}\underset{g}{}\ldots\underset{g}{}b_{k}\underset{g}{}(\mathcal{D}%
_{a}^{+}X). \label{MCD.5iv}%
\end{align}

The verification of Eq.(\ref{MCD.5iii}) is trivial.

To verify Eq.(\ref{MCD.5iv}) we will use \emph{complete induction }over the
$k$ smooth vector fields $b_{1},\ldots,b_{k}.$

Let us take $b\in\mathcal{V}(U),$ by using Eq.(\ref{MCD.5ii}) and
Eq.(\ref{MCD.5i}), we have
\begin{align}
\mathcal{D}_{a}^{+}(b\underset{g}{}X)  &  =\mathcal{D}_{a}^{+}(b\underset
{g}{\lrcorner}X)+\mathcal{D}_{a}^{+}(b\wedge X)\nonumber\\
&  =(\mathcal{D}_{a}^{+}b)\underset{g}{\lrcorner}X+b\underset{g}{\lrcorner
}(\mathcal{D}_{a}^{+}X)+(\mathcal{D}_{a}^{+}b)\wedge X+b\wedge(\mathcal{D}%
_{a}^{+}X),\nonumber\\
\mathcal{D}_{a}^{+}(b\underset{g}{}X)  &  =(\mathcal{D}_{a}^{+}b)\underset
{g}{}X+b\underset{g}{}(\mathcal{D}_{a}^{+}X). \label{MCD.5v}%
\end{align}

Now, let us $b_{1},\ldots,b_{k},b_{k+1}\in\mathcal{V}(U).$ By using twice the
inductive hypothesis and Eq.(\ref{MCD.5v}) we can write
\begin{align*}
&  \mathcal{D}_{a}^{+}(b_{1}\underset{g}{}\ldots\underset{g}{}b_{k}%
\underset{g}{}b_{k+1}\underset{g}{}X)\\
&  =(\mathcal{D}_{a}^{+}(b_{1}\underset{g}{}\ldots\underset{g}{}%
b_{k}))\underset{g}{}b_{k+1}\underset{g}{}X+b_{1}\underset{g}{}\ldots
\underset{g}{}b_{k}\underset{g}{}(\mathcal{D}_{a}^{+}(b_{k+1}\underset{g}%
{}X))\\
&  =(\mathcal{D}_{a}^{+}(b_{1}\underset{g}{}\ldots\underset{g}{}%
b_{k}))\underset{g}{}b_{k+1}\underset{g}{}X+b_{1}\underset{g}{}\ldots
\underset{g}{}b_{k}\underset{g}{}(\mathcal{D}_{a}^{+}b_{k+1})\underset{g}{}X\\
&  +b_{1}\underset{g}{}\ldots\underset{g}{}b_{k}\underset{g}{}b_{k+1}%
\underset{g}{}(\mathcal{D}_{a}^{+}X)\\
&  =((\mathcal{D}_{a}^{+}(b_{1}\underset{g}{}\ldots\underset{g}{}%
b_{k}))\underset{g}{}b_{k+1}+b_{1}\underset{g}{}\ldots\underset{g}{}%
b_{k}\underset{g}{}(\mathcal{D}_{a}^{+}b_{k+1}))\underset{g}{}X\\
&  +b_{1}\underset{g}{}\ldots\underset{g}{}b_{k}\underset{g}{}b_{k+1}%
\underset{g}{}(\mathcal{D}_{a}^{+}X)\\
&  =(\mathcal{D}_{a}^{+}(b_{1}\underset{g}{}\ldots\underset{g}{}b_{k}%
\underset{g}{}b_{k+1}))\underset{g}{}X+b_{1}\underset{g}{}\ldots\underset{g}%
{}b_{k}\underset{g}{}b_{k+1}\underset{g}{}(\mathcal{D}_{a}^{+}X).\blacksquare
\end{align*}

From the theory of extensors (\cite{3},\cite{4}) it follows that given a
metric field $g,$ there is a \emph{non-singular} smooth $(1,1)$-extensor field
$h$ such that
\begin{equation}
g=h^{\dagger}\circ\eta\circ h, \label{MCD.6}%
\end{equation}
where $\eta$ is an \emph{orthogonal metric field} with the same signature as
$g.$ Such $h$ will be said to be a \emph{gauge metric field} for $g.$

It might as well be asked whether there is any relationship between a
$g$-compatible pair of $a$-\emph{DCDO's} and a $\eta$-compatible pair of
$a$-\emph{DCDO's.} The answer is YES.

\textbf{Theorem 2. }Let $h$ be a gauge metric field for $g.$ For any
$g$-compatible pair of $a$-\emph{DCDO's,} namely $(_{g}\mathcal{D}_{a}^{+}$
$,_{g}\mathcal{D}_{a}^{-}),$ there exists an unique $\eta$-compatible pair of
$a$-\emph{DCDO's,} namely $(_{\eta}\mathcal{D}_{a}^{+},_{\eta}\mathcal{D}%
_{a}^{-}),$ such that\footnote{Recall that $h^{*}=(h^{-1})^{\dagger
}=(h^{\dagger})^{-1},$ and $\underline{h}^{-1}=(\underline{h})^{-1}%
=\underline{(h^{-1})}$ and $\underline{h}^{\dagger}=(\underline{h})^{\dagger
}=\underline{(h^{\dagger})}.$}
\begin{align}
\underline{h}(_{g}\mathcal{D}_{a}^{+}X)  &  =\text{ }_{\eta}\mathcal{D}%
_{a}^{+}\underline{h}(X),\label{MCD.7a}\\
\underline{h}^{\ast}(_{g}\mathcal{D}_{a}^{-}X)  &  =\text{ }_{\eta}%
\mathcal{D}_{a}^{-}\underline{h}^{\ast}(X). \label{MCD.7b}%
\end{align}
And reciprocally, given any $\eta$-compatible pair of $a$-\emph{DCDO's}, say
$(_{\eta}\mathcal{D}_{a}^{+},_{\eta}\mathcal{D}_{a}^{-}),$ there is an unique
$g$-compatible pair of $a$-\emph{DCDO's,} say $(_{g}\mathcal{D}_{a}^{+}%
,_{g}\mathcal{D}_{a}^{-}),$ such that the above formulas are satisfied.

\textbf{Proof}

Given $(_{g}\mathcal{D}_{a}^{+},_{g}\mathcal{D}_{a}^{-}),$ since $h$ is a
non-singular smooth $(1,1)$-extensor field, we can indeed construct a
well-defined pair of $a$-\emph{DCDO's,} namely $(_{h}\mathcal{D}_{a}^{+}%
,_{h}\mathcal{D}_{a}^{-}),$ by the following formulas
\[
_{h}\mathcal{D}_{a}^{+}X=\underline{h}(_{g}\mathcal{D}_{a}^{+}\underline
{h}^{-1}(X))\text{ and }_{h}\mathcal{D}_{a}^{-}X=\underline{h}^{\ast}%
(_{h}\mathcal{D}_{a}^{+}\underline{h}^{\dagger}(X)).
\]
As defined above, $(_{h}\mathcal{D}_{a}^{+},_{h}\mathcal{D}_{a}^{-})$ is the
$h$-\emph{deformation} of $(_{g}\mathcal{D}_{a}^{+},_{g}\mathcal{D}_{a}^{-})$.

But, it is obvious that $_{h}\mathcal{D}_{a}^{+}$ and $_{h}\mathcal{D}_{a}%
^{-}$ as defined above satisfy in fact Eq.(\ref{MCD.7a}) and Eq.(\ref{MCD.7b}%
), i.e., $\underline{h}(_{g}\mathcal{D}_{a}^{+}X)=$ $_{h}\mathcal{D}_{a}%
^{+}\underline{h}(X)$ and $\underline{h}^{\ast}(_{g}\mathcal{D}_{a}^{-}X)=$
$_{\eta}\mathcal{D}_{a}^{-}\underline{h}^{\ast}(X).$

In order to check the $\eta$-compatibility of $(_{h}\mathcal{D}_{a}^{+}%
,_{h}\mathcal{D}_{a}^{-}),$ we can write
\begin{align*}
(_{h}\mathcal{D}_{a}^{++}\eta)(b)  &  =\text{ }_{h}\mathcal{D}_{a}^{-}%
\eta(b)-\eta(_{h}\mathcal{D}_{a}^{+}b)\\
&  =h^{\ast}(_{g}\mathcal{D}_{a}^{-}h^{\dagger}\circ\eta(b))-\eta\circ
h(_{g}\mathcal{D}_{a}^{+}h^{-1}(b))\\
&  =h^{\ast}(_{g}\mathcal{D}_{a}^{-}h^{\dagger}\circ\eta\circ h\circ
h^{-1}(b)-h^{\dagger}\circ\eta\circ h(_{g}\mathcal{D}_{a}^{+}h^{-1}(b)))\\
&  =h^{\ast}(_{g}\mathcal{D}_{a}^{-}g(h^{-1}(b))-g(_{g}\mathcal{D}_{a}%
^{+}h^{-1}(b))),\\
&  =h^{\ast}(_{g}\mathcal{D}_{a}^{++}g)(h^{-1}(b)).
\end{align*}
This implies that $_{h}\mathcal{D}_{a}^{++}\eta=h^{\ast}\circ(_{g}%
\mathcal{D}_{a}^{++}g)\circ h^{-1}$. Then, since $_{g}\mathcal{D}_{a}%
^{++}g=0,$ it follows that $_{h}\mathcal{D}_{a}^{++}\eta=0,$ i.e.,
$(_{h}\mathcal{D}_{a}^{+},_{h}\mathcal{D}_{a}^{-})$ is $\eta$-compatible.

Now, if there exists another $\eta$-compatible pair $(_{\eta}\mathcal{D}%
_{a}^{\prime+},_{\eta}\mathcal{D}_{a}^{\prime-})$, which satisfies
Eq.(\ref{MCD.7a}) and Eq.(\ref{MCD.7b}), i.e.,
\[
\underline{h}(_{g}\mathcal{D}_{a}^{+}X)=\text{ }_{\eta}\mathcal{D}_{a}%
^{\prime+}\underline{h}(X)\text{ and }\underline{h}^{\ast}(_{g}\mathcal{D}%
_{a}^{-}X)=\text{ }_{\eta}\mathcal{D}_{a}^{\prime-}\underline{h}^{\ast}(X).
\]
Then, by substituting $X$ for $\underline{h}^{-1}(X)$ in the first one and $X$
for $\underline{h}^{\dagger}(X)$ in the second one, it follows that $_{\eta
}\mathcal{D}_{a}^{\prime+}=$ $_{\eta}\mathcal{D}_{a}^{+}$ and $_{\eta
}\mathcal{D}_{a}^{\prime-}=$ $_{\eta}\mathcal{D}_{a}^{-}$.

So the existence and uniqueness are proved. Such a $\eta$-compatible pair of
$a$-\emph{DCDO's} satisfying Eq.(\ref{MCD.7a}) and Eq.(\ref{MCD.7b}) is just
the $h$-deformation of the $g$-compatible pair of $a$-\emph{DCDO's.}

By following analogous steps we can also prove that such a $g$-compatible pair
of $a$-\emph{DCDO's} satisfying Eq.(\ref{MCD.7a}) and Eq.(\ref{MCD.7b}) is
just the $h^{-1}$-deformation of the $\eta$-compatible pair of $a$%
-\emph{DCDO's}.$\blacksquare$

In the eight paper of this series we study the relation between the curvature
and torsion tensors of a pair $(_{\eta}\mathcal{D}_{a}^{+},_{\eta}%
\mathcal{D}_{a}^{-})$ and its deformation $(_{g}\mathcal{D}_{a}^{+}%
,_{g}\mathcal{D}_{a}^{-})$.

\section{Conclusions}

In this paper the geometric calculus on manifolds developed in previous papers
\cite{1,2} is used to introduce through the concept of a metric extensor field
$g$ a metric structure on a smooth manifold $M$. The associated Levi-Civita
connection extensor \textit{field} and Christoffel \textit{operators} and are
given. The concept of a geometrical structure for a manifold $M$ as a triple
$(M,g,\gamma)$, where $\gamma$ is a connection extensor field defining a
parallelism structure for $M$ and the associated theory of metric compatible
covariant derivatives is also introduced. Eventually the most important
results of the present paper are: (i): the decomposition of the Levi-Civita
operator of the first kind (Eq.(\ref{LCC.1})) and, (ii): the determination of
a relationship between the connection extensor fields and covariant
derivatives corresponding to deformed geometrical structures (Theorem 2).
Those results combined with the results of the seventh and eight papers on
this series are important in the formulation and applications of \ geometrical
theories of the gravitational field, as we will show elsewhere.

\textbf{Acknowledgments: }V. V. Fern\'{a}ndez and A. M. Moya are very grateful
to Mrs. Rosa I. Fern\'{a}ndez who gave to them material and spiritual support
at the starting time of their research work. This paper could not have been
written without her inestimable help.

\end{document}